\documentclass[12pt]{article}
\usepackage{amssymb,amsmath, amsthm, amscd,ifthen}
\usepackage[T2A]{fontenc}
\usepackage[cp1251]{inputenc}
\usepackage[english,russian]{babel}
\usepackage[dvips]{graphicx}
\usepackage[hyper]{amsbib}
\usepackage{indentfirst}
\usepackage{amsmath,amssymb}
\usepackage{amsmath,amsthm,amsfonts,amssymb,latexsym,enumerate,graphicx,parskip}
\usepackage{marginnote, etoolbox, imakeidx, needspace}

\theoremstyle{plain}
\newtheorem{theorem}{Theorem}

\newtheorem{lemma}{Lemma}
\newtheorem{claim}{Claim}
\newtheorem{corollary}{Corollary}

\newenvironment{proof1}{\noindent{\it Proof.\,}}{\hfill$\Box$}

 \date{} 
\begin{document}
	
\title{Chain method for panchromatic colorings of hypergraphs}

\author{Margarita Akhmejanova\footnote{Moscow Institute of Physics and Technology, Laboratory of Combinatorial and Geometric Structures, Laboratory of Advanced Combinatorics and Network Applications, 141700, Institutskiy per. 9, Dolgoprudny, Moscow Region, Russia. E-mail: mechmathrita@gmail.com;}\label{1}, \quad J\'ozsef Balogh\footnote{Department of Mathematical Sciences, University of Illinois at Urbana-Champaign, Urbana, Illinois
61801, USA. E-mail: jobal@illinois.edu;}}

\maketitle
\textbf{Abstract.} We deal with an extremal problem concerning  panchromatic colorings of hypergraphs. A vertex $r$-coloring of a hypergraph $H$ is \emph{panchromatic} if every edge meets every color. We prove that for every $r<\sqrt[3]{\frac{n}{100\ln n}}$, every $n$-uniform hypergraph $H$ with $|E(H)|\leq \frac{1}{20r^2}\left(\frac{n}{\ln n}\right)^{\frac {r-1}{r}}\left(\frac{r}{r-1}\right)^{n-1}$ has a panchromatic coloring with $r$ colors. 

\section{Introduction and related work}

We study colorings of uniform hypergraphs. Let us recall some definitions.

\emph{A vertex $r$-coloring} of a hypergraph $H=(V,E)$ is a mapping from the vertex set $V$ to a set of $r$ colors. An $r$-coloring of $H$ is \emph{panchromatic} if  each edge has at least one vertex of each color. 

The first sufficient condition on the existence of a panchromatic coloring of a hypergraph was obtained in 1975 by Erd\H{o}s and Lov\'asz \cite{ErdLov}.
They proved that if every edge of an $n$-uniform hypergraph intersects at
most  
\begin{equation}\label{erdlov}
   \frac{r^{n-1}}{4(r-1)^n}
\end{equation}
other edges then the hypergraph  has a panchromatic
coloring with $r$ colors. 

The next generalization of the problem was formulated in 2002 by Kostochka \cite{Kost}, who posed the following question:\emph{
What is the minimum possible number of edges in an $n$-uniform hypergraph that does not admit a panchromatic coloring with $r$ colors?} He denoted this number by $p(n,r)$. 

Following closely behind this problem is a related one: a hypergraph $H=(V,E)$ has property $B$ if there is a coloring of $V$ by $2$ colors so that no edge $f \in E$ is monochromatic. Erd\H{o}s and Hajnal \cite{Erd2} (1961) proposed  to find the value $m(n)$ equal to the minimum possible number of edges in a $n$-uniform hypergraph without property $B$. Erd\H{o}s \cite{Erd} (1963--1964) found  bounds $\Omega\left(2^{n}\right) \leq m(n)=O\left(2^{n} n^{2}\right)$ and Radhakrishnan and Srinivasan \cite{RadhSrin} $(2000)$ proved  $m(n)\geq\Omega\left(2^{n}(n / \ln n)^{1 / 2}\right)$.
Clearly, $m(n)=p(n,2)$.

We return to the panchromatic coloring. Kostochka \cite{Kost} has found connections between $p(n,r)$ and  minimum possible number of vertices in a $k$-partite graph with list chromatic number greater than $r$. Using results of Erd\H{o}s, Rubin and Taylor \cite{ErdRubTey} and also Alon's result \cite{Alon} Kostochka \cite{Kost} proved the existence of constants $c_1$ and $c_2$ that for every  large $n$ and fixed $r$:
\begin{equation}\label{label_Kost}
\frac{e^{c_{1}\frac{n}{r}}}{r}\leq p(n,r)\leq r e^{c_{2}\frac{n}{r}}.
\end{equation}

In 2010, bounds (\ref{label_Kost}) were considerably improved in the paper of Shabanov \cite{Shab4}:
$$
p(n, r) \geqslant \frac{\sqrt{21}-3}{4 r}\left(\frac{n}{(r-1)^{2} \ln n}\right)^{1 / 3}\left(\frac{r}{r-1}\right)^{n},~~~~\text{for all~~$r<n,$}
$$
$$
p(n, r) \leqslant \frac{1}{r}\left(\frac{r}{r-1}\right)^{n} e(\ln r) \frac{n^{2}}{2(r-1)} \varphi_{1},~~~~\text{when ~~$r=o(\sqrt{n}),$}
$$
$$
p(n, r) \leqslant \frac{1}{r}\left(\frac{r}{r-1}\right)^{n} e(\ln r) n^{3 / 2} \varphi_{2},~~~~\text{when~~$n=o\left(r^{2}\right),$}
$$
where $\varphi_{1},\varphi_{2}$ some functions of $n$ and $r(n)$, tending to one at $n\rightarrow\infty$.

\bigskip
In 2012, Rozovskaya and Shabanov \cite{RosShab} improved Shabanov's lower bound by proving that for $r<n /(2 \ln n)$
\begin{equation}\label{bound:RosShab}
\frac{1}{2r^2}\left(\frac{n}{\ln n}\right)^{1 / 2}\left(\frac{r}{r-1}\right)^{n} \leqslant p(n, r) \leqslant c_{2} n^{2}\left(\frac{r}{r-1}\right)^{n} \ln r.
\end{equation}
Further research was conducted by Cherkashin \cite{Cherk} in 2018. In his work, Cherkashin introduced the auxiliary value $p^{'}(n, r)$, which is numerically equal to the minimum number of edges in the class of $n$-uniform hypergraphs
$H = (V, E)$, in which any subset of vertices $V^{\prime} \subset V \text {with}\left|V^{\prime}\right| \geq\left[\frac{r-1}{r}|V|\right]$ must contain an edge. Analyzing the value $p^{'}(n, r)$ and using Sidorenko's \cite{Sid} estimates on the Turan numbers, Cherkashin proved that for $n\geq 2, r\geq 2$
$$
p(n, r) \leq c \frac{n^{2} \ln r}{r}\left(\frac{r}{r-1}\right)^{n}.
$$
Cherkashin also proved that for $r \leq c \frac{n}{\ln n}$
\begin{equation}\label{Check}
p(n, r) \geq c \max \left(\frac{n^{1 / 4}}{r \sqrt{r}}, \frac{1}{\sqrt{n}}\right)\left(\frac{r}{r-1}\right)^{n}.
\end{equation}

And repeating the ideas of Gebauer \cite{Geb} Cherkashin constructed an example of a hypergraph that has few edges and does not admit a panchromatic coloring in $r$ colors. The reader is referred to the survey \cite{surv} for the detailed history of panchromatic colorings.

It is thus natural to consider the local case. Formally, \emph{the degree of an edge  $A$} is the number of hyperedges intersecting $A$. Let $d(n, r)$ be the minimum possible value of the maximum edge degree in an $n$-uniform hypergraph that does not admit panchromatic coloring with $r$ colors. Then, the Erd\H{o}s and Lov\'asz  result (\ref{erdlov}) can be easily translated into following form:
\begin{equation}\label{erdlov2}
d(n,r)\geq\frac{r^{n-1}}{4(r-1)^n}.
\end{equation}
However, the bound (\ref{erdlov2}) appeared not to be sharp. The restriction on $d(n, r)$ have been improved by Rozovskaya and  Shabanov \cite{RosShab}. In their work they achieved that
\begin{equation}\label{Roz_Shab_local_bound}
	d(n,r)>\frac{\sqrt{11}-3}{4 r(r-1)}\left(\frac{n}{\ln n}\right)^{1 / 2}\left(\frac{r}{r-1}\right)^{n},~~~~~\text{when $r \leqslant n /(2 \ln n).$}
\end{equation}

\section{Our results}

The main result of our paper improves the estimate \eqref{bound:RosShab} as follows.

\bigskip
\begin{theorem}\label{thm:main}
Suppose $r\leq\sqrt[3]{\frac{n}{100\ln n}}$. Then we have
\begin{equation}
p(n,r)\geq \frac{1}{20r^2}\left(\frac{n}{\ln n}\right)^{\frac {r-1}{r}}\left(\frac{r}{r-1}\right)^{n}.
\end{equation}
\end{theorem}

\bigskip
\begin{corollary}
	There is an absolute constant $C$ so that for every $n>2$ and  $\ln n<r<\sqrt[3]{\frac{n}{100\ln n}}$
	\begin{equation*}
	p(n,r)\geq \frac{Cn}{r^2 \ln n}\cdot e^{\frac{n}{r}+\frac{n}{2r^2}}.
	\end{equation*}
\end{corollary}

We refine the bound \eqref{Roz_Shab_local_bound} as follows.

\bigskip
\begin{theorem}\label{thm:main2}
	For every $2<r<\sqrt[3]{\frac{n}{100\ln n}}$ 
	\begin{equation}
	d(n,r)\geq\frac{1}{40r^3}\left(\frac{n}{\ln n}\right)^{\frac {r-1}{r}}\left(\frac{r}{r-1}\right)^{n}.
	\end{equation}
\end{theorem}

\subsection{Methods}

In the work, we propose a new idea based on the Pluhar ordered chain method \cite{Pluhar}.
In the case of panchromatic coloring, the resulting structure is no longer a real ordered chain, but rather an intricate "snake ball". Nevertheless, with the help of probabilistic analysis, we managed to obtain a strong lower bound.

The rest of this paper is organised as follows. The next section describes a coloring algorithm. Section \ref{Section 4} is devoted to the detailed analysis of the algorithm. In Section \ref{Section 5} we collect some inequalities that will be subsequently useful. The last two sections contain proofs of Theorems~\ref{thm:main} and \ref{thm:main2}.

\section{The coloring algorithm}

We may and will assume that $r \geq 3$, because case $r=2$ corresponds to the case $m(n)$. Let $H=(V, E)$ be an $n$-uniform hypergraph with less than $\frac{1}{20r^2}\left(\frac{n}{\ln n}\right)^{\frac {r-1}{r}}\left(\frac{r}{r-1}\right)^{n}$ edges and let $r<\sqrt[3]{\frac{n}{100\ln n}}$. We will show that $H$ has a panchromatic coloring with $r$ colors. 

\begin{figure}[h]
	\includegraphics[width=1.\linewidth]{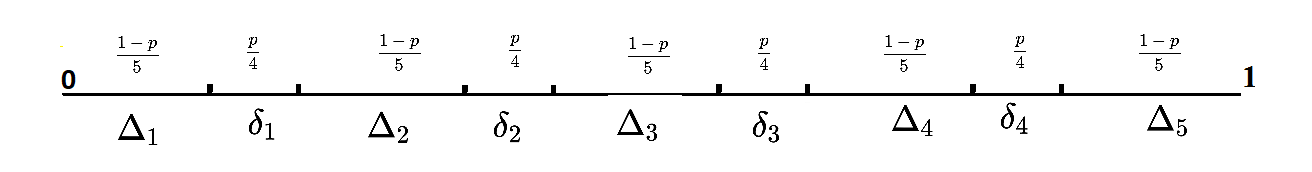}
	\selectlanguage{english}
\caption{Partition of $[0,1)$ into $\Delta_1,\delta_1,\Delta_2,\delta_2,\ldots,\Delta_5$ when $r=5$.}
\label{picture18}
\end{figure}

\bigskip
We define a special random order on the set $V$ of vertices of hypergraph $H$ using a mapping $\sigma: V \rightarrow[0,1],$ where $\sigma(v), v \in V$ -- i.i.d. with uniform distribution on $[0,1]$. The value $\sigma(v)$ we will call the \emph{weight} of the vertex $v$. Reorder the vertices so that $\sigma(v_{1})<\ldots<\sigma(v_{|V|})$. Put
\begin{equation}\label{choice_p}
p=\left(\frac{r-1}{r}\right)\frac{(r-1)^2\ln(\frac{n}{\ln n})}{n}.
\end{equation}
We divide the unit interval $[0,1)$ into subintervals $\Delta_1,\delta_1,\Delta_2,\delta_2,\ldots,\Delta_r$ as on the Figure \ref{picture18}, i.e.
$$
  \Delta_i=\left[(i-1)\left(\frac {1-p}r+\frac p{r-1}\right),i\cdot\frac {1-p}r+(i-1)\cdot\frac p{r-1}\right),\;i=1,\ldots,r;
$$
$$
  \delta_i=\left[i\cdot\frac {1-p}r+(i-1)\cdot\frac p{r-1},i\left(\frac {1-p}r+\frac p{r-1}\right)\right),\;i=1,\ldots,r-1.
$$
The length of each large subinterval $\Delta_i$ is equal to $\frac{1-p}{r}$ and every small subinterval $\delta_i$ has length equal to $\frac{p}{r-1}$. Since $p<\frac{1}{100 r}$ under the given assumptions on $r$, we can see that the intervals $\Delta_{1}, \ldots, \Delta_{r}$ are each wider than the intervals $\delta_{1} \ldots, \delta_{r-1}$. A vertex $v$ is said to belong to a subinterval $[c,d)$, if $\sigma(v)\in[c,d)$. We note that the same division of the segment $[0,1]$ has already been used by the first author for proving some bounds on proper colorings \cite{AkhShab}.

\bigskip

\bigskip
We color the vertices of hypergraph $H$ according to the following algorithm, which consists of two steps.
\begin{enumerate}
\item  First, each $v\in\Delta_{i}$ is colored with color $i$ for every $i\in[r].$
\item  Then, moving with the growth of $\sigma$, we color a vertex $v\in\delta_{i}$ with color $i$ if there exists an edge $e, v\in e$ such that $e$ does not have color $i$ in the current coloring. Otherwise we color $v$ with color $i+1$.
\end{enumerate}

\section{Analysis of the algorithm}\label{Section 4}
\subsection{Short edge}\label{Section 4.1}

We say that an edge $A$ is \emph{short} if $A\cap (\Delta_i\cup\delta_{i})=\emptyset $ or $A\cap (\Delta_{i+1}\cup\delta_{i})=\emptyset$ for some $i\in[r-1]$. The probability of this event for fixed edge $A$ and fixed $i$ is at most $2\left(1-\left(\frac{1-p}{r}+\frac{p}{r-1}\right)\right)^n$. Summing up this upper bound over all edges and $i \in[r-1]^{n}$ we get 

\begin{gather*}\label{prob_1}
2(r-1)|E|\left(1-\left(\frac{1-p}{r}+\frac{p}{r-1}\right)\right)^n\leq\frac{2(r-1)}{20r^2}\left(\frac{n}{\ln n}\right)^{\frac {r-1}{r}}\left(\frac{r}{r-1}\right)^n\cdot\\
\cdot\left(\frac{r-1}{r}-\frac{p}{r(r-1)}\right)^n
\leq\frac{1}{10r}\left(\frac{n}{\ln n}\right)^{\frac {r-1}{r}}\left(1-\frac{p}{(r-1)^2}\right)^n\leq\frac{1}{10r}.
\end{gather*}

Hence, we conclude that the expected number of short edges is less than $1/10r$, hence with probability at least $1-1/10r$ there is no short edge. 

\subsection{Snake ball} \label{Section 4.2}

Suppose our algorithm fails to produce a panchromatic $r$-coloring and there is no short edges. Let $A$ be an edge, which does not contain some color $i$. 

Now we have two possibilities:
\begin{itemize}
\item $i<r,$ in this situation edge $A$ is disjoint from the interval $\Delta_i\cup\delta_i$, which means that $A$ is short, a contradiction.
\item $i=r.$ 

\begin{figure}[h]
	\includegraphics[width=1.\linewidth]{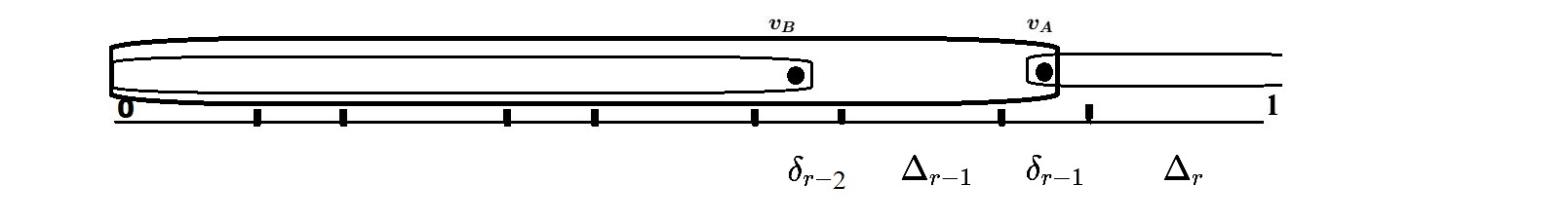}
	\selectlanguage{english}
    \caption{Edges $A$ and $B$ in a snake ball.}
    \label{picture}
\end{figure}
\end{itemize}

Edge $A$ is not short, so $A\cap(\delta_{r-1}\cup\Delta_r)\ne\emptyset$. Since $A$ does not contain color $r$ we have $A\cap\Delta_r=\emptyset$. Denote $v_A$ the last vertex of $A\cap\delta_{r-1}$. We note that $v_A$ could receive color $r-1$ only if at the moment of coloring $v_A$ there was an edge $B$ without color $r-1$ and $v_A$ was the first vertex of $B\cap\delta_{r-1}$. In this situation we say that the pair $(A,B)$ is \emph{conflicting in $\delta_{r-1}$} and the vertex $v_A$ is \emph{dangerous vertex in $\delta_{r-1}$.}

Again, edge $B$ is not short and did not contain color $r-1$ at the moment of coloring $v_A$, so $B\cap(\delta_{r-2}\cup\Delta_{r-1})\ne\emptyset$ and $B\cap\Delta_{r-1}=\emptyset$. For $v_B$, the last vertex of  $B\cap\delta_{r-2}$ , there exists an edge $C$, which at the moment of coloring $v_B$ was without color $r-2$ and $v_B$ was the first vertex of $C\cap\delta_{r-2}$. We get $(B,C)$ is conflicting pair in $\delta_{r-2}$ and $v_B$ is \emph{dangerous vertex in $\delta_{r-2}$.}

\bigskip
Repeating the above arguments, we obtain a construction called \emph{snake ball}. It is an edge sequence $H'=(C_1=A,C_2=B,...,C_{r})$ such that consecutive edges $(C_i,C_{i+1})$ form conflicting pairs in $\delta_{r-i}$. 

\bigskip
Summarizing the above, we can say that

\bigskip
\begin{claim}
 If for injective $\sigma:V\to[0;1)$ there are neither snack balls nor short edges then Algorithm 1 produces a panchromatic $r$-coloring.
\end{claim}

\bigskip
\begin{lemma}\label{lemma_prob1}
Let $H'=(C_1,\ldots,C_r)$ be an ordered $r$-tuple of edges in the hypergraph $H$.
Then the probability of the event that $H'$ forms a snake ball and all the edges $C_1,\ldots,C_r$ are not short does not exceed
$$
\left(\frac{p}{r-1}\right)^{r-1}\left(\frac{r-1}{r}\right)^{(n-2)r}\prod_{v\in H':s(v)\geq 2
}\frac{\left(1-s(v)\frac{1-p}{r}\right)}{\left(1-\left(\frac{1-p}{r}+\frac{2p}{r-1}\right)\right)^{s(v)}}\prod_{i=1}^{r-1}|C_{i}\cap C_{i+1}|,
$$
where $s(v)$ is the number of edges of $H'$ that contain vertex $v$.
\end{lemma}

Before we present the proof of this lemma, we introduce some facts and  give the basic scheme of the proof. Note that if $v\in C_i$ then $\sigma(v)\notin\Delta_{r-i+1}$. Furthermore, for each $v$ its weight $\sigma(v)$  belongs to the subintervals of total length at most
\begin{equation}\label{property_s(v)}
1-s(v)\frac{1-p}{r}.
\end{equation}
The scheme of the proof is following:
\begin{itemize}
\item fix vertex $v_j\in C_j\cap C_{j+1}$ and its weight $\sigma(v_j)$ for all $j=1,\ldots, r-1$. Assuming that $v_j$ is the dangerous vertex in $\delta_{r-j}$ calculate conditional probability given weights of dangerous vertices.
\item sum up (integrate) the previous probability over all possible values of weights, using that $\sigma(v_{j})\in \delta_{r-j}$, as this is needed for $H^{\prime}$ to be a snake ball.
\item Finally, sum over all choices of $v_{1}, \ldots, v_{r-1}$.
\end{itemize}

\begin{proof1} 
Fix dangerous vertex $v_j\in C_j\cap C_{j+1}$ for each $j=1,\ldots r-1.$ 
Put $[\alpha_j,\beta_j)=\delta_j$, $\beta_j-\alpha_j=p/(r-1)$ and $y_j=\beta_{r-j}-\sigma(v_j)$. Recall that $0 \leq y_{j} \leq p /(r-1).$

Fix for a moment variables $y_1,\ldots,y_{r-1}$. Then, for $v\in C_i$ with $s(v)=1$ its weight $\sigma(v)$ belongs to the subinterval of total length at most
$$
1-\left(\frac{1-p}{r}+y_{i+1}+\frac{p}{r-1}-y_i\right)~~~\text{if}~~ i\in[2,r-1].
$$ 
And similarly, $1-\left(\frac{1-p}{r}+y_1\right)$ for $i=1$ and $1-\left(\frac{1-p}{r}+\frac{p}{r-1}-y_{r-1}\right)$
for $i=r$.

Now we are ready to give an upper bound for the probability of the event that ``$H'$ forms a snake ball'', conditional on the value taken by $y_1,\ldots,y_{r-1}$:

\begin{gather}
\left(1-\left(\frac{1-p}{r}+y_1\right)\right)^{n-1}\cdot\left(1-\left(\frac{1-p}{r}+y_2+\frac{p}{r-1}-y_1\right)\right)^{n-2}\cdot\ldots\cdot\\
\cdot\left(1-\left(\frac{1-p}{r}+y_{r-1}+\frac{p}{r-1}-y_{r-2}\right)\right)^{n-2}\cdot\left(1-\left(\frac{1-p}{r}+\frac{p}{r-1}-y_{r-1}\right)\right)^{n-1}\cdot\\
\cdot\prod_{v\in H':s(v)\geq 2}\frac{\left(1-s(v)\frac{1-p}{r}\right)}{\left(1-\left(\frac{1-p}{r}+\frac{2p}{r-1}\right)\right)^{s(v)}}.
\end{gather}

\bigskip
Here we estimated as if all the rest of the vertices have $s(v)=1$ (factors (11) and factor (12)), and then using (\ref{property_s(v)}), edited for vertices with $s(v)>1$ by multiplying by $1-s(v)\frac{1-p}{r}$ and divided by $\left(1-\left(\frac{1-p}{r}+\frac{2p}{r-1}\right)\right)^{s(v)}.$ The factor $\left(1-\left(\frac{1-p}{r}+\frac{2p}{r-1}\right)\right)$ is obviously no more than any factor for $s(v)=1$, so we get a correct upper bound.

Taking out factor $\left((r-1)/r\right)^{(n-2)r+2}$ in the above equation and using estimate $(1+y)^s\leq\exp\{ys\}$, we get the following upper bound on product of (11) and (12):
\begin{align*}
&\left(\frac{r-1}{r}\right)^{r(n-2)+2} \exp \left(\frac{(n-1) p}{r-1}-\frac{(n-2) p}{r-1}-\frac{p}{(r-1)^2}-\frac{r y_{1}}{r-1}+\frac{r y_{r-1}}{r-1}\right)\leq \\
&\left(\frac{r-1}{r}\right)^{r(n-2)+2} \exp \left(\frac{p(r-2)}{(r-1)^2}+\frac{r y_{r-1}}{r-1}\right)\leq\left(\frac{r-1}{r}\right)^{r(n-2)+2} \exp \left(\frac{p(r-2)}{(r-1)^2}+\frac{r p}{(r-1)^{2}}\right)=\\
&\left(\frac{r-1}{r}\right)^{r(n-2)+2} \exp \left(\frac{2p}{r-1}\right)<\left(\frac{r-1}{r}\right)^{r(n-2)}.
\end{align*}
To obtain the final estimate, we have to integrate over the weights $y_1,y_2,\ldots,y_{r-1}$ (factor $\left(p/(r-1)\right)^{r-1}$) and
sum up over all possible choices for the  $v_{1},...,v_{r-1}$ (factor $\prod_{i=1}^{r-1}|C_i\cap C_{i+1}|$).

\end{proof1}

\section{Auxilary calculations}\label{Section 5}

Under the assumptions of
Theorem \ref{thm:main} we will formulate and prove three auxiliary lemmas needed to prove Theorem \ref{thm:main}. In particular, in Lemma 2, we replace product of pairwise intersections on their sum $\sum_{i<j}|C_i\cap C_j|$ and in Lemma 4, we will use double-counting for estimating the sum $\sum_{i<j}|C_i\cap C_j|$, which can be large with $n$, by special bounded terms.

\bigskip
\begin{lemma}\label{auxilary_lemma_1}
Let $H'=(C_1,\ldots,C_r)$ be an ordered $r$-tuple of edges in the hypergraph $H$. Then 
\begin{equation}
\sum_{\pi \in S_{r}}|C_{i_1}\cap C_{i_2}||C_{i_2}\cap C_{i_3}|\cdot \ldots\cdot|C_{i_{r-1}}\cap C_{i_r}|\leq\left(\frac{2\sum_{i<j}|C_i\cap C_j|+r}{r}\right)^r,
\end{equation}
where $S_r$ denotes all permutations $\pi=(i_1,\ldots,i_r)$ of $(1,2,\ldots,r)$.
\end{lemma}

\begin{proof1} Denote the cardinality of the edge intersection $|C_i\cap C_j|$ by $x_{i,j}$. Then, we have to prove that

$$
\sum_{\pi \in S_{r}}x_{i_1,i_2}x_{i_2,i_3}\cdot\ldots\cdot x_{i_{r-1},i_r}
\leq\left(\frac{2\sum_{i<j} x_{i,j}+r}{r}\right)^r.
$$

First, we will show that
\begin{align}\label{permutation}
&\sum_{\pi \in S_{r}}x_{i1,i2}x_{i_2,i_3}\cdot \ldots\cdot x_{i_{r-1},i_r}
\leq(x_{1,2}+\ldots+x_{1,r}+1)\cdot\ldots\cdot(x_{r,1}+\ldots+x_{r,r-1}+1).
\end{align}

Let us call $(x_{i,1}+\ldots+x_{i,r}+1)$ from (\ref{permutation}) the \emph{bracket number $i$}. We define a mapping $f$ between elements from the left-hand side of (\ref{permutation}) and ordered sets that are obtained after performing the multiplication in (\ref{permutation}).

Let $f:x_{i1,i2}x_{i_2,i_3}\ldots x_{i_{r-1},i_r}\mapsto x_{1,t_1}x_{2,t_2}\ldots x_{r,t_r}$, where $x_{1,t_1}x_{2,t_2}\ldots x_{r,t_r}$ is the product of the following $r$ elements: $x_{i_{r-1},i_r}$ from the bracket number $i_{r-1}$,  $x_{i_{r-2},i_{r-1}}$ from the bracket number $i_{r-2}$ and so forth, finally we take the factor 1 from the unused bracket. For example,

 $x_{5,6}x_{6,1}x_{1,4}x_{4,3}x_{3,2}$ is mapped to $x_{1,4}\cdot 1\cdot x_{3,2}\cdot x_{4,3}\cdot x_{5,6}\cdot x_{6,1}$. 

 We note that $f$ is an injection. Indeed, 
for each $x_{1,t_1}x_{2,t_2}\ldots x_{r,t_r}$ there exists at most one sequence $x_{i1,i2}x_{i_2,i_3}\ldots x_{i_{r-1},i_r}$, with $i_1\neq i_2\ldots\neq i_r$, such as $f(x_{i1,i2}x_{i_2,i_3}\ldots x_{i_{r-1},i_r})=x_{1,t_1}\ldots x_{r,t_r}$.

 So, since $f$ does not change the product and $f$ is an injection we get that the right-hand side of (\ref{permutation}) is not less than the left-hand side.

Finally, by the inequality on the arithmetic-geometric means and by $x_{i,j} =  x_{j,i}$

$$(x_{1,2}+\ldots+x_{1,r}+1)\cdot\ldots\cdot(x_{r,1}+\ldots+x_{r,r-1}+1)\leq\left(\frac{2\sum_{i<j} x_{i,j}+r}{r}\right)^r.$$
\end{proof1}

\begin{lemma}\label{auxilary_lemma_2}
For all $s\in\{2,\ldots,r-1\}$
\begin{equation}
\frac{\left(1-s\frac{1-p}{r}\right)}{\left(1-\left(\frac{1-p}{r}+\frac{2p}{r-1}\right)\right)^{s}}\leq e^{-\frac{s^2}{20r^2}}.
\end{equation}
\end{lemma}

\begin{proof1}
	First prove the case $s\geq 3$. 
	\begin{gather}\label{lemma_2_exp2}
	\frac{\left(1-\frac{s(1-p)}{r}\right)}{\left(1-\left(\frac{1-p}{r}+\frac{2p}{r-1}\right)\right)^{s}}=\frac{\left(1-\frac{s(1-p)}{r}\right)}{\left(1-\left(\frac{1-p}{r}\right)\right)^{s}\left(1-\frac{2pr}{(r-1)(r-1+p)}\right)^s}\leq\frac{\left(1-\frac{s(1-p)}{r}\right)\left(1+\frac{1-p}{r-1+p}\right)^s}{ \left(1-\frac{2pr}{(r-1)^2}\right)^s}.
	\end{gather}
Now we deal with factors in (\ref{lemma_2_exp2}) separetely:
\begin{gather*}
\left(1+\frac{1-p}{r-1+p}\right)^s\leq \left(1+\frac{1-p}{r-1}\right)^s=
\text{|Apply Taylor's formula with Lagrange Remainder|=}\\
1+\frac{s(1-p)}{r-1}+\frac{s(s-1)(1-p)^2}{2(r-1)^2}+
\frac{s(s-1)(s-2)(1-p)^3(1+\theta\cdot\frac{1-p}{r-1})^{s-3}}{6(r-1)^3}\leq\\
\text{bound $(s-1)/(r-1)$ by $s/r$, $(s-1)(s-2)/(r-1)^2$ by$s^2/r^2$ and $(1+\theta/(r-1))^{s-3}$ by $e$.}\\
\leq 1+\frac{s(1-p)}{r-1}+\frac{s^2(1-p)}{2r(r-1)}+
\frac{s^3(1-p)^2 e}{6r^2(r-1)}.
\end{gather*}
 Hence, the numerator of (\ref{lemma_2_exp2}) does not exceed
\begin{gather*}
\left(1-\frac{s(1-p)}{r}\right)\left(1+\frac{s(1-p)}{r-1}+\frac{s^2(1-p)}{2r(r-1)}+
\frac{s^3(1-p)^2}{2r^2(r-1)}\right)<
1-\frac{s^2(1-p)}{r(r-1)}
\left(1-p-1/2\right)+\\\frac{s(1-p)}{r(r-1)}=1-\frac{s^2(1-p)}{r(r-1)}\left(1/2-1/s-p\right)<1-\frac{s^2(1/6-p)(1-p)}{r^2}<1-\frac{s^2}{7r^2}\leq\exp\left\{-\frac{s^2}{7r^2}\right\}.
\end{gather*}
Using bounds $1/(1-x)<1+2x$ for $x<1/2$ and estimating $pr<1/100$, which follows from restrictions on $r$, we finally get
	\begin{gather*}
\frac{\left(1-s\frac{1-p}{r}\right)}{\left(1-\left(\frac{1-p}{r}+\frac{2pr}{r-1}\right)\right)^{s}}\leq\exp\left\{-\frac{s^2}{7r^2}\right\}\left(1-\frac{2pr}{(r-1)^2}\right)^{-s}<\exp\left\{-\frac{s^2}{7r^2}\right\}\left(1+\frac{4pr}{(r-1)^2}\right)^s\leq\\
\exp\left\{\frac{4prs}{(r-1)^2}-\frac{s^2}{7r^2}\right\}\leq\exp\left\{\frac{s}{25(r-1)^2}-\frac{s^2}{7r^2}\right\}<\exp\left\{\frac{4}{25s}\cdot\frac{s^2}{r^2}-\frac{s^2}{7r^2}\right\}<\exp\left\{-\frac{s^2}{20r^2}\right\}.
	\end{gather*}
Consider the case $s=2$.
\begin{gather*}
	\frac{1-2(1-p)/r}{\left(1-\left(\frac{1-p}{r}+\frac{2p}{r-1}\right)\right)^{2}}\leq\frac{1-2(1-p)/r}{1-\frac{2(1-p)}{r}-\frac{4p}{(r-1)}+\frac{1}{2r^2}}\leq \frac{1-2(1-p)/r}{1-\frac{2(1-p)}{r}-\frac{1}{4r^2}+\frac{1}{2r^2}}=1-\frac{1/4r^2}{1-2\frac{1-p}{r}+\frac{1}{4r^2}}\\
	\leq 1-1/4r^2\leq\exp\{-1/4r^2\}<\exp\{-1/5r^2\},
\end{gather*}
where we used that $4p/(r-1)<8p/r=8pr/r^2<8/100r^2<1/4r^2.$
\end{proof1}

\bigskip
\begin{lemma}\label{auxilary_lemma_3}
	\begin{gather}\label{statement_lemma_4}
	\left(\prod_{v\in H':s(v)\geq 2}\frac{\left(1-s(v)\frac{1-p}{r}\right)}{\left(1-\left(\frac{1-p}{r}+\frac{2p}{r-1}\right)\right)^{s(v)}}\right)\cdot\sum_{\sigma \in S_{r}}|C_{i_1}\cap C_{i_2}||C_{i_2}\cap C_{i_3}|\cdot\ldots\cdot |C_{i_{r-1}}\cap C_{i_r}|\leq 20^rr^{2r}e^{-r+1}
	\end{gather}
\end{lemma}

\begin{proof1}
By Lemmas \ref{auxilary_lemma_1} and \ref{auxilary_lemma_2} the left hand side of (\ref{statement_lemma_4}) does not exceed
\begin{align*}
&\exp\left\{-\sum_{v\in H':s(v)\geq 2}\frac{s^2(v)}{20r^2}\right\}\left(\frac{2\sum_{i<j}|C_i\cap C_j|+r}{r}\right)^r.
\end{align*}
Now we will use the following double-counting: $\sum_{i<j}|C_i\cap C_j|$ is equal to $\sum_{v\in H':s(v)\geq 2} \binom{s(v)}{2}<1/2\sum_{v\in H':s(v)\geq 2}s^2(v)$. Hence,
\begin{align*}
&\exp\left\{-\sum_{v\in H':s(v)\geq 2}\frac{s^2(v)}{20r^2}\right\}\left(\frac{2\sum_{i<j}|C_i\cap C_j|+r}{r}\right)^r\leq\exp\left\{-\sum_{v\in H':s(v)\geq 2}\frac{s^2(v)}{20r^2}\right\}\cdot r^r\cdot\\
&\cdot\left(\frac{\sum_{v\in H':s(v)\geq 2} s^2(v)+r}{r^2}\right)^r\leq r^re^{-t/20}(t+1)^r\leq\frac{20^rr^{2r}}{e^{r-1}},
\end{align*}
where we used $t=\sum_{v\in H':s(v)\geq 2} s^2(v)/r^2$ and observed that 
the expression $\left((t+1)^{r} e^{-t / 20}\right)$ is maximized when $t=20r-1$.
\end{proof1}

\section{Proof of Theorem \ref{thm:main}}\label{section6}

We want to show that there is a positive probability that no edge is short and no tuple of edges form a snake ball.

\bigskip
Denote $\sum^{*}$ the sum over all $r$-sets $J\subseteq (1,2,\ldots,|E|)$, $\sum^{o}$ 
the sum over all ordered $r$-tuples $(j_1,\ldots,j_r)$, with $\{j_1,\ldots,j_r\}$ forming such a $J$ and $\sum_{\pi\in S_r}$ denote the sum over all permutations $\pi=(i_1,\ldots,i_r)$ of $(1,2,\ldots,r)$.

In Section \ref{Section 4.1} we already proved that the expected number of short edges does not exceed $1/(10r)$. The expected number of snake ball can be upper bounded as follows: 
\begin{align*}
&\stackrel{o}{\sum}\mathbb{P}\left(\left(C_{j_1},...,C_{j_r}\right)\text{~forms a snake ball}\right)
=\stackrel{*}{\sum}\sum_{\pi\in S_r}\mathbb{P}\left(\left(C_{i_1},...,C_{i_r}\right)\text{~forms a snake ball}\right).\\
&\text{On the other hand,}\\
&\sum_{\pi\in S_r}\mathbb{P}\left(\left(C_{i_1},...,C_{i_r}\right)\text{~forms a snake ball }\right)\\
&\leq\sum_{\pi\in S_r}\left(\frac{p}{r-1}\right)^{r-1}\left(\frac{r-1}{r}\right)^{(n-2)r}\prod_{v\in H':s(v)\geq 2}\frac{\left(1-s(v)\frac{1-p}{r}\right)}{\left(1-\left(\frac{1-p}{r}+\frac{2p}{r-1}\right)\right)^{s(v)}}|C_{i_1}\cap C_{i_{2}}|\ldots|C_{i_{r-1}}\cap C_{i_{r}}|\\
&=\left(\frac{p}{r-1}\right)^{r-1}\left(\frac{r-1}{r}\right)^{(n-2)r}\prod_{v\in H':s(v)\geq 2}\frac{\left(1-s(v)\frac{1-p}{r}\right)}{\left(1-\left(\frac{1-p}{r}+\frac{2p}{r-1}\right)\right)^{s(v)}}\sum_{\pi\in S_r}|C_{i_1}\cap C_{i_{2}}|\ldots|C_{i_{r-1}}\cap C_{i_{r}}|\\
&\leq\left(\frac{p}{r-1}\right)^{r-1}\left(\frac{r-1}{r}\right)^{(n-2)r}\frac{20^rr^{2r}}{e^{r-1}}\leq\left(\frac{(r-1)^2\ln(\frac{n}{\ln n})}{rn}\right)^{r-1}\cdot\left(\frac{r-1}{r}\right)^{(n-2)r}\cdot\frac{20^rr^{2r}}{e^{r-1}},
\end{align*}
where for the first inequality we used Lemma \ref{lemma_prob1} and for the second Lemma \ref{auxilary_lemma_3} and in the final inequality we took $p$ from \ref{choice_p}. Finally,
\begin{align*}
&\stackrel{*}{\sum}\sum_{\pi\in S_r}\mathbb{P}\left(\left(C_{i_1},...,C_{i_r}\right)\text{~forms a snake ball}\right)\leq\\
&\binom{|E|}{r}\cdot\left(\frac{(r-1)^2\ln(\frac{n}{\ln n})}{rn}\right)^{r-1}\cdot\left(\frac{r-1}{r}\right)^{(n-2)r}\cdot\frac{20^rr^{2r}}{e^{r-1}}\leq\\
&\frac{\left(\frac{1}{20r^2}\left(\frac{n}{\ln n}\right)^{\frac {r-1}{r}}\left(\frac{r}{r-1}\right)^{n}\right)^r}{r!}\cdot\left(\frac{(r-1)^2\ln(\frac{n}{\ln n})}{rn}\right)^{r-1}\cdot\left(\frac{r-1}{r}\right)^{(n-2)r}\cdot\frac{20^rr^{2r}}{e^{r-1}}\leq\frac{1}{r}\left(\frac{r}{r-1}\right)^2.
\end{align*}

Since $1-\frac{1}{10r}-\frac{1}{r}\left(\frac{r}{r-1}\right)^2>0$, with positive probability the Algorithm creates a panchromatic coloring with $r$ colors, which proves Theorem \ref{thm:main}.

\bigskip
\begin{corollary}
	There is an absolute constant $c$ so that for every $n>2$ and  $\ln n<r<\sqrt[3]{\frac{n}{100\ln n}}$
	\begin{equation*}
	p(n,r)\geq c\frac{n}{r^2(\ln n)}\cdot e^{\frac{n}{r}+\frac{n}{2r^2}}.
	\end{equation*}
\end{corollary}

\begin{proof1}
By applying Taylor's formula with Peano remainder, we obtain 
$$
\left(1+\frac{1}{r-1}\right)e^{-\frac{1}{r}-\frac{1}{2r^2}}=1+\frac{1}{3r^3}+O\left(\frac{1}{r^4}\right).
$$
Thus, $\left(1+\frac{1}{r-1}\right)>e^{\frac{1}{r}+\frac{1}{2r^2}}$. Finally, we use $\left(\frac{n}{\ln n}\right)^{-\frac{1}{r}}>\frac{1}{e}$ when $r>\ln n$ and Theorem \ref{thm:main}.
\end{proof1}

\section{Local variant: proof of Theorem \ref{thm:main2}}\label{section 7}

A useful parameter of $H$ is its \emph{maximal edge degree} $$D: =D(H) = \max_{e \in E(H)}\left|\left\{e^{\prime} \in E(H): e \cap e^{\prime} \neq 0\right\}\right|.$$
We show that for $3<r<\sqrt[3]{\frac{n}{100\ln n}}$ every $n$-uniform hypergraph with $D\leq\frac{1}{40r^3}\left(\frac{n}{\ln n}\right)^{\frac {r-1}{r}}\left(\frac{r}{r-1}\right)^{n}$ has a panchromatic coloring with $r$ colors, which implies Theorem \ref{thm:main2}.

\bigskip
Let us recall Lov\'asz Local Lemma, which shows a
useful sufficient condition for simultaneously avoiding a set $A_1, A_2,\ldots,A_N$ of ``bad'' events:

\bigskip
\begin{lemma}[The Local Lemma; General Case, \cite{ErdLov}]\label{LLL}
Let $A_{1}, A_{2}, \ldots, A_{n}$ be events in an arbitrary probability space. A directed graph $\overline{D}=(V, E)$ on the set of vertices $V=\{1,2, \ldots, n\}$ is  $a$ dependency digraph for the events $A_{1}, \ldots, A_{n}$ if for each
i, $1 \leq i \leq n$, the event $A_{i}$ is mutually independent of all the events $\left\{A_{j}:(i, j) \notin E\right\}$. Suppose that $\overline{D}=(V, E)$ is a dependency digraph for the above events and suppose there are real numbers $x_{1}, \ldots, x_{n}$ such that $0 \leq x_{i}<1$ and $\mathbb{P}\left[A_{i}\right] \leq x_{i} \prod_{(i, j) \in E}\left(1-x_{j}\right)$
for all $1 \leq i \leq n .$ Then
\[
\mathbb{P}\left[\bigwedge_{i=1}^{n} \overline{A_{i}}\right] \geq \prod_{i=1}^{n}\left(1-x_i\right).
\]
In particular, with positive probability, no event $A_{i}$ holds.
\end{lemma}

 To prove Theorem \ref{thm:main2} we will use the following generalization of Lemma \ref{LLL}.

\bigskip
\begin{lemma}\label{LLL_special_case}
If all events have probability $\mathbb{P}\left(A_{i}\right) \leq \frac{1}{2}$, and for all $i$
\begin{equation}\label{LLL_special_case_1/4}
\sum_{j:(i,j)\in E}\mathbb{P}\left(A_{j}\right) \leq \frac{1}{4},
\end{equation}
 then there is a positive probability that no $A_{i}$ holds.
\end{lemma}

For the sake of completeness, we give the proof of Lemma \ref{LLL_special_case} here.

\begin{proof1}
Put $x_{i}=2\mathbb{P}\left(A_{i}\right).$ Then, for all $i$ 
\[x_{i} \prod_{(i, j) \in E}\left(1-x_{j}\right)=2\mathbb{P}\left(A_{i}\right) \prod_{(i, j) \in E}\left(1-2\mathbb{P}(A_j)\right)\geq \mathbb{P}\left(A_{i}\right).\]
\end{proof1}

In our case the set of bad events has two types: short edges and snake balls. Let $\mathcal{Q}(C)$ be the event ``edge $C$ is short'' and $\mathcal{W}(C_1,\ldots,C_r)$ be the event ``$(C_1,\ldots,C_r)$ forms a snake ball and all the edges $C_1,\ldots,C_r$ are not short''. Note that $\mathcal{Q}(C)$ depends on at most on $D+1$ events $\mathcal{Q}(C')$ and at most on $2r(D+1)D^{r-1}$ events $\mathcal{W}(C_1,\ldots,C_r)$. Similarly, $\mathcal{W}(C_1,\ldots,C_r)$ depends at most on $r(D+1)$ events $\mathcal{Q}(C')$ and at most on $2r^2(D+1)D^{r-1}$ events $\mathcal{W}(C'_1,\ldots,C'_r)$. Hence, using bounds from Sections \ref{Section 4.1} and \ref{section6} we get the following upper bounds:
\begin{enumerate}
\item if $A_i=\mathcal{W}(C_1,\ldots,C_r):$
\begin{align*}
&\sum_{j:(i,j)\in E}\mathbb{P}(A_j)\leq r(D+1)\cdot 2(r-1)\left(1-\left(\frac{1-p}{r}+\frac{p}{r-1}\right)\right)^n+\\
+2r^2&(D+1) D^{r-1}\cdot\left(\frac{r-1}{r}\right)^{(n-2)r}\left(\frac{p}{r-1}\right)^{r-1}\frac{20^rr^{2r}}{e^{r-1}}<\frac{2r^2}{40r^3}+\frac{2r^2}{r2^re^{r-1}}<\frac{1}{4}.
\end{align*}
\item if $A_i=\mathcal{Q}(C):$
\begin{align*}
&\sum_{j:(i,j)\in E}\mathbb{P}(A_j)\leq (D+1)\cdot 2(r-1)\left(1-\left(\frac{1-p}{r}+\frac{p}{r-1}\right)\right)^n+\\
&+2r(D+1) D^{r-1}\cdot\left(\frac{r-1}{r}\right)^{(n-2)r}\left(\frac{p}{r-1}\right)^{r-1}\frac{20^rr^{2r}}{e^{r-1}}<\frac{1}{4}.
\end{align*}
\end{enumerate}
In both cases inequality (\ref{LLL_special_case_1/4}) holds, completing the proof of Theorem \ref{thm:main2}.

\section{Acknowledgements}

The work of the first  was funded by RFBR, project number 20-31-70039 (Theorem \ref{thm:main}) and the Council for the Support of
Leading Scientific Schools of the President of the Russian Federation, grant no. N.Sh.-2540.2020.1 (Theorem \ref{thm:main2}). The first author is a Young Russian Mathematics award winner and would like to thank its
sponsors and jury.

\renewcommand{\refname}{References}


\begin{thebibliography}{99}\label{bibliography}

\bibitem{AkhShab} M. Akhmejanova, D.A. Shabanov, ``Equitable colorings of hypergraphs with few edges'', \emph{Discrete Applied Mathematics}, \textbf{276}, 2020, 2--12.

\bibitem{Alon} N. Alon, ``Choice number of graphs: a probabilistic approach'', \emph{Combin. Probab. Comput.},
\textbf{1}:2, 1992, 107--114.

\bibitem{Cherk} D. Cherkashin, ``A note on panchromatic colorings'', \emph{Discrete Mathematics}, \textbf{341}:3, 2018, 652--657. 

\bibitem{surv} Raigorodskii, A., Cherkashin, D, ``Extremal problems in hypergraph colourings'', \emph{ Russian Mathematical Surveys}, \textbf{75}:1, 2020, 89--146. 

\bibitem{CherkKozik} D. Cherkashin, J. Kozik, ``A note on random greedy coloring of uniform hypergraphs'', \emph{Random Struct. Alg.}, \textbf{47}:3, 2015, 407--413.

\bibitem{Erd} P. Erd\H{o}s. ``On a combinatorial problem, II. In J. Spencer, editor, Paul Erd\H{o}s: The Art of Counting'', \emph{MIT Press}, 1973, 445--447.

\bibitem{Erd2}  P. Erd\H{o}s, A. Hajnal, ``On a property of families of sets'', \emph{Acta Math. Hung.}, \textbf{12}:1, 1961, 87--123.

\bibitem{ErdLov} P. Erd\H{o}s, L. Lov\'asz, ``Problems and results on 3-chromatic hypergraphs and some related questions'', \textit{Infinite and Finite Sets}, Coll Math Soc J Bolyai, \textbf{10}, 1975, 609--627.

\bibitem{ErdRubTey} P. Erd\H{o}s, A. Rubin, H. Taylor, ''Choosability in graphs'', \textit{Proceedings of the West
Coast Conference on Combinatorics, Graph Theory and Computing (Humboldt State Univ., Arcata, Calif., 1979)} , Congr. Numer., \textbf{26}, Utilitas Math. Publ., Winnipeg,
Man., 1980, 125--157.

\bibitem{Geb} H. Gebauer, ``On the construction of 3-chromatic hypergraphs with few edges'', \emph{Journal of Combinatorial Theory, Series A}, \textbf{120}:7, 2013, 1483--1490. 

\bibitem{Kost} A. Kostochka, ``On a theorem of Erd\H{o}s, Rubin, and Taylor on choosability of complete bipartite graphs'', \emph{Electron. J. Combin.}, \textbf{9}:1, 2002, 1--4.

\bibitem{Pluhar} A. Pluh\'ar, ``Greedy colorings for uniform hypergraphs'', \emph{ Random Struct. Alg.}, \textbf{35}:2, 2009, 216--221.

\bibitem{RadhSrin} J. Radhakrishnan, A. Srinivasan, ``Improved bounds and algorithms for hypergraph two-coloring'', \emph{Random Struct. Alg.}, \textbf{16}:1, 2000, 4--32.

\bibitem{RosShab} A. Rozovskaya, D. Shabanov, ``Extremal problems for panchromatic colourings of uniform hypergraphs'', \emph{Discrete Math. Appl.}, \textbf{22}:2, 2012, 185--206.

\bibitem{Shab4} D. Shabanov, ``On a generalization of Rubin's theorem.'',   \emph{J. Graph Theory}, \textbf{67}:3, 2011, 226--234.

\bibitem{Sid} A. Sidorenko, ``What we know and what we do not know about Turan numbers'', \emph{Graphs and Combinatorics}, \textbf{11}:2, 1995, 179--199.

\end{thebibliography}
\end{document}